\documentclass[11pt,a4paper]{article}
\usepackage[T1]{fontenc}
\usepackage{lmodern}
\usepackage{authblk}
\usepackage{graphicx}
\usepackage{amsmath,amssymb,amsthm,mathtools}
\usepackage{booktabs,tabularx,enumitem}
\usepackage[margin=0.89in]{geometry}
\usepackage[final]{microtype}
\usepackage[colorlinks=true,linkcolor=blue,citecolor=blue,urlcolor=blue]{hyperref}
\newtheorem{proposition}{Proposition}
\newtheorem{corollary}{Corollary}
\newtheorem{remark}{Remark}
\newcommand{\F}{F_\mu}
\newcommand{\G}{G_\mu}
\newcommand{\Ps}{\Psi_\mu}
\newcommand{\Exp}{\operatorname{exp}}
\newcommand{\Log}{\operatorname{log}}
\newcommand{\Cov}{\operatorname{Cov}}
\newcommand{\iid}{\mathrel{\overset{\mathrm{i.i.d.}}{\sim}}}
\title{\vspace{-1.9em}\textbf{A Note on the Cut-Locus-Corrected CLT for local empirical Fr\'echet Means }}
\author{Susovan Pal}
\affil{Visiting Research Fellow, Hausdorff Research Institute of Mathematics, University of Bonn\\ \texttt{susovan97@gmail.com}}
\date{\today}
\begin{document}
\maketitle
\vspace{-1.6em}
\begin{abstract}
We record three observations concerning \cite{HLW24}. First, the proof of its Theorem~1 yields the population linearization (13)--(14) at \textit{any critical point} of the Fr\'echet function $F_\mu$, under its remaining geometric and measure-theoretic assumptions; neither \textit{isolation} nor a \textit{nonsingular} Hessian of $F_\mu$ is needed. Second, its Theorem~2 is given a shorter proof by applying the M-estimator theorem \cite[Theorem~5.23, pp.~53--54]{vdV98} to negative squared distance in normal coordinates, as opposed to the Donsker argument used there. Third, the above shorter route also helps us establish a CLT for a measurable selection of strict local empirical Fréchet mean with positive definite Hessian, with the cut locus correction, relaxing the assumptions of uniqueness or global minimality, as in \cite{HLW24}. Numerical experiments at the end show the agreement of the Hessian of the Fréchet function at \textit{local} Fréchet means with the (signed) sum of the expected Hessian of the squared distance function and the cut locus correction there.
\end{abstract}
\vspace{-0.55em}
\section{Setting and correspondence of notation}
Let $(M,g)$ be a compact connected $m$-dimensional Riemannian manifold without boundary, with distance $\rho$, and let $\xi_1,\xi_2,\ldots\iid\mu$. Write
\[
\F(x)=\int_M\rho^2(x,y)\,d\mu(y),\qquad
F_n(x)=\frac1n\sum_{i=1}^n\rho^2(x,\xi_i),\qquad
\G(x)=\int_M\Log_x(y)\,\mathbf 1_{\{y\notin C_x\}}\,d\mu(y).
\]
The first and third quantities are (1) and (3) of \cite[pp.~1221--1222]{HLW24}; $F_n$ corresponds to (2), and $C_x$ denotes the cut locus. Under assumptions (i)--(v) of \cite[Theorem~1, p.~1224]{HLW24}, with $x_0$ its \textit{unique, global} Fr\'echet mean, the parallel-transport expansion (13)--(14) holds with the \textit{total Hessian} operator $\Ps(x_0)$, including the cut-locus correction $J_\mu(x_0)$. If in addition $\Ps(x_0)$ is positive definite, \cite[Theorem~2, p.~1225]{HLW24} states, for any measurable global empirical minimizer $\hat x_n$,
\begin{equation}\label{eq:hlwclt}
\sqrt n\,\Log_{x_0}(\hat x_n)\xrightarrow{\mathcal D}
N_m\!\left(0,\Ps(x_0)^{-1}V_0\Ps(x_0)^{-\mathsf T}\right),\quad
V_0=\Cov(\Log_{x_0}\xi_1).
\end{equation}
Our use of \cite[Theorem~5.23]{vdV98} relies on the following dictionary that help us translate the general loss function framework there into the Fréchet function framework in \cite{HLW24}. Fix an orthonormal basis of $T_pM$ and a sufficiently small open Euclidean ball $\Theta=B_r(0)\subset T_{p}M$ on which $\Exp_{p}$ is a diffeomorphism onto its image. In \cite{HLW24}, $p$ is the unique, global Fréchet mean, but we will relax this to be a critical point of $F_\mu$ in our Proposition \ref{prop:stat}.
\vspace{-0.35em}
\begin{center}
\small
\begin{tabularx}{\linewidth}{@{}l>{\raggedright\arraybackslash}X>{\raggedright\arraybackslash}X@{}}
\toprule
\cite{vdV98}, Theorem 5.23 & Geometric specialization & Meaning \\
\midrule
$\mathcal X$, $X_i$, $P$, $\mathbb P_n$ & $M$, $\xi_i$, $\mu$, $n^{-1}\sum_i\delta_{\xi_i}$ & data space, observations, laws \\
$\Theta\subset\mathbb R^m$, $\theta_0$ & $B_r(0)\subset T_pM$, $0$ & (open) descriptor space, optimizer \\
$m_\theta(y)$ & $-\rho^2(\Exp_{p}\theta,y)$ & empirical criterion/loss \\
$Pm_\theta$, $\mathbb P_nm_\theta$ & $-\F(\Exp_{p}\theta)$, $-F_n(\Exp_{p}\theta)$ & population, empirical objectives/losses \\
$\hat\theta_n$ & $\Log_{p}\hat x_n$ (eventually) & estimator in normal coordinates \\
$\dot m_{\theta_0}(y)$ & $2\Log_{p} y$ ($\mu$-a.e.) & criterion derivative/$L^2$-envelope \\
$V_{\theta_0}$ & $-2\Ps({p})$ & population Hessian \\
\bottomrule
\end{tabularx}
\end{center}
Note that, the parameter $\theta$ determines a \emph{candidate location}; the observation $y$ remains in \textit{all of $M$}. In particular, choosing $\Theta$ does \emph{not} exclude observations near $C_{p}$, see Remark \ref{rem:choice-of-theta-0}.
%%%%%%%%%%%%%%%%%%%%%%%%%%%%%%%%%%%%%%%%%%%%%%%%%%%%%%%%%%%%%%%%%%%%%%%%%%%%%%%%%%%%%%%%%%%%%%%%%%%%%%%%%%%%%%%%%%%%%%%%%%%%%%%%
\section{Generalization of Theorem 1 of \cite{HLW24} at stationary points}

In this section, we show that the linearization formula (13) and the cut-locus corrected Hessian formula (14) from \cite{HLW24} remains valid for \textit{any critical point} of $F_\mu.$

\begin{proposition}\label{prop:stat}
Let $p\in M$ satisfy assumptions (i), (iii), (iv), and (v) of \cite[Theorem~1, p.~1224]{HLW24}, with $p$ in place of $x_0$. If $\nabla\F(p)=0$, then, as $x\to p$,
\begin{equation}\label{eq:lin}
\Pi_{x,p}\G(x)=-\Ps(p)\Log_p(x)+o(\rho(x,p)),\qquad
\Ps(p)=-\int_M H(p\mid y)\,d\mu(y)-J_\mu(p),
\end{equation}
where $\Pi_{x,p}$ is parallel transport from $T_xM$ to $T_pM$; the $(1,1)$-operator $H$ given by negative one half of the Hessian of the squared distance from $p,$ and the cut locus correction $J_\mu$ are respectively defined by (8) and (20) of \cite{HLW24}. Note that $p$ is only a critical point: there is no assumption that $p$ is a local/global minimum, isolated, or that $\Ps(p)$ is nonsingular.
\end{proposition}
\begin{proof}
The geometric Proposition~3 and the estimates in the proof of Theorem~1 of \cite[pp.~1231--1237, especially p.~1234 after (25)]{HLW24} use assumption (ii) only to conclude $\G(x_0)=0$. Under assumption (iii), the cut loci have $\mu$-measure zero near $p$, and the first-variation identity gives $\nabla\F(p)=-2\G(p)$; hence stationarity supplies exactly that step.\footnote{We use the factor $\nabla F_\mu=-2G_\mu$ dictated by (1), (3), and (5) of \cite{HLW24}; the displayed proportionality immediately preceding (15) on p.~1225 of that article has the numerical factor reversed.} The remainder of the cited proof, including the calculation of $J_\mu$ in (26)--(28), is unchanged. No inverse of $\Ps(p)$ occurs. Thus (13)--(14) of \cite{HLW24} hold with $p$ replacing $x_0$.
\end{proof}
Along radial geodesics, \eqref{eq:lin} and $\nabla\F=-2\G$ yield the \emph{pointwise} second-order expansion
\begin{equation}\label{eq:peano}
\F(\Exp_p v)=\F(p)+\langle\Ps(p)v,v\rangle+o(\|v\|^2),\qquad v\to0.
\end{equation}
Consequently $\operatorname{Hess}^M\F(p)[u,v]=2\langle\Ps(p)u,v\rangle$ for $u,v\in T_pM$, as in (15) of \cite[p.~1225]{HLW24}. No $C^2$ assumption in a neighborhood is needed. In normal coordinates the Riemannian and coordinate Hessians agree \emph{at the center}; no interchange of the Hessian with the integral defining $\F$ is used.

\section{A direct M-estimation proof of Theorem~2 in \cite{HLW24} and its extension to \textit{local} empirical Fr\'echet means with positive Hessians}\label{sec:direct-M-estimation-proof}

We first verify the hypotheses of \cite[Theorem~5.23, pp.~53--54]{vdV98} with $p=x_0$, thereby recovering \cite[Theorem~2]{HLW24} for global empirical Fr\'echet means. We then extend this result to a measurable selection of local empirical Fr\'echet minimizers converging to a strict local minimum $p$ of $\F$, under the hypotheses of Proposition~\ref{prop:stat} and the positive definiteness of $\Ps(p)$, \textit{without requiring uniqueness or global minimality of $p$}; see Corollary~\ref{cor:local}.

\begin{enumerate}[label=\textup{(\alph*)},leftmargin=*,itemsep=2pt,topsep=3pt]
\item \textbf{Differentiability almost surely.} Assumption (iii) of \cite[Theorem~1]{HLW24} implies $\mu(C_p)=0$. For $y\notin C_p$, first variation of squared distance gives
\[D_\theta m_\theta(y)|_{\theta=0}[v]=2\langle\Log_p y,v\rangle,\quad \dot m_0(y)=2\Log_p y.\]
By compactness these derivatives have finite second moments, and their expectation is zero because $p$ is stationary.
\item \textbf{Local Lipschitz square-integrable envelope.} Put $D=\operatorname{diam}M$. On a smaller convex Euclidean parameter ball, $\Exp_p$ has a finite Lipschitz constant $L$, so for all $y\in M$ and $\theta,\eta\in\Theta$,
\[
|m_\theta(y)-m_\eta(y)|\le 2D\,\rho(\Exp_p\theta,\Exp_p\eta)
\le 2DL\|\theta-\eta\|.
\]
Thus the constant envelope $2DL$ lies in $L^2(\mu)$; no support restriction is used.
\item \textbf{Quadratic population criterion with negative definite Hessian.} From \eqref{eq:peano},
\[
Pm_\theta=-\F(\Exp_p\theta)=Pm_0+\tfrac12\theta^{\mathsf T}V_{\theta_0}\theta+o(\|\theta\|^2),\qquad V_{\theta_0}=-2\Ps(p).
\]
Assumption (vi) of \cite[Theorem~2]{HLW24} makes this matrix negative definite. Crucially, $\Ps(p)$ is the \textit{cut-locus corrected,  full Hessian} in (14) in \cite{HLW24}, including $J_\mu(p)$.
\item \textbf{Consistency and maximization.} Proposition~1 of \cite[p.~1222]{HLW24} yields $\hat x_n\to p$ almost surely for any measurable global empirical minimizer. Hence eventually $\hat x_n\in\Exp_p(\Theta)$, and $\hat\theta_n=\Log_p\hat x_n\to0$. On this event $\hat\theta_n$ \emph{exactly} maximizes $\mathbb P_nm_\theta$ over $\Theta$. Defining $\hat\theta_n=0$ on the exceptional event gives a sequence in $\Theta$ satisfying the $o_P(n^{-1})$ near-maximization condition (for the precise definition, see \cite{vdV98}, Theorem 5.23, pp. 53 with $o_P(n^{-1})$, and not the \textit{weaker} condition on pp. 45 with $o_P(1)$ bound); the exceptional event has probability tending to zero, and thus $\Pr(n\Delta_n>\varepsilon)\le\Pr(\hat x_n\notin\Exp_p\Theta)\to0$ for every $\varepsilon>0$, 
where $\Delta_n:=\sup_{\theta\in\Theta}P_nm_\theta-P_nm_{\hat\theta_n}\ge0$.

\end{enumerate}
Theorem~5.23 therefore gives
\begin{align*}
\sqrt n\,\hat\theta_n
&=-V_{\theta_0}^{-1}\frac1{\sqrt n}\sum_{i=1}^n\dot m_0(\xi_i)+o_P(1)
=\Ps(p)^{-1}\frac1{\sqrt n}\sum_{i=1}^n\Log_p(\xi_i)+o_P(1).
\end{align*}
Since $\Cov(\dot m_0(\xi_1))=4V_0$, the ordinary multivariate CLT proves \eqref{eq:hlwclt}. Restricting the parameter chart does not alter the integral over $M$ or remove its cut-locus correction.

\begin{remark}[\textbf{On the choice of $\Theta_0$}]\label{rem:choice-of-theta-0}
    For the setting of the unique global Fréchet mean in \cite{HLW24}, set $p=x_0$, the unique minimizer of $\F$ on $M$. Since $p\notin C_p$, the same point uniquely minimizes $\F$ over $M\setminus C_p$. The tangential non-cut domain $\Theta_{\max}:=\Log_p(M\setminus C_p)$ is open, and $\Exp_p:\Theta_{\max}\to M\setminus C_p$ is a diffeomorphism. Hence $0=\operatorname*{argmin}_{\theta\in\Theta_{\max}}\F(\Exp_p\theta)$, or equivalently $0=\operatorname*{argmax}_{\theta\in\Theta_{\max}}Pm_\theta$, where $Pm_\theta=-\int_M\rho^2(\Exp_p\theta,y)\,d\mu(y)$. The maximal domain is unnecessary for \cite[Theorem~5.23, p.~53]{vdV98}: any fixed, sufficiently small ball $\Theta=B_r(0)\subset\Theta_{\max}$ works because $0$ remains a population maximizer. Likewise, the empirical estimator is originally a global minimizer $\hat x_n\in\operatorname*{argmin}_{q\in M}F_n(q)$. By \cite[Proposition~1, p.~1222]{HLW24}, $\hat x_n\to p$ almost surely; hence eventually a.s. $\hat x_n\in\Exp_p(\Theta)$, and $\hat\theta_n=\Log_p\hat x_n$ exactly maximizes $\mathbb P_nm_\theta=-F_n(\Exp_p\theta)$ over $\Theta$. On the exceptional event $\hat x_n\notin \Exp_p(\Theta)$  we set $\hat\theta_n=0$, as in part~(d), without affecting the $o_P(n^{-1})$ near-maximization condition. Thus $\Theta$ restricts the candidate locations, not the observations $y\in M$, and does not remove the cut-locus correction in $\Ps(p)$.
\end{remark}

\begin{corollary}[\textbf{CLT for \textit{local} empirical Fréchet means}]\label{cor:local}
Suppose $p$ satisfies the hypotheses of Proposition~\ref{prop:stat} and $\Ps(p)$ is positive definite. For a sufficiently small $r>0$, let $\hat p_n$ be a measurable global minimizer of $F_n$ over the closed geodesic ball $\overline B_r(p)$. Then $\hat p_n\to p$ a.s.; eventually $\hat p_n$ is an interior, hence local, empirical Fr\'echet minimizer; and we have the cut-locus corrected CLT for $\hat p_n$:
\[
\sqrt n\,\Log_p(\hat p_n)\xrightarrow{\mathcal D}
N_m\!\left(0,\Ps(p)^{-1}\Cov(\Log_p\xi_1)\Ps(p)^{-\mathsf T}\right).
\]
\end{corollary}

\begin{proof}
By \eqref{eq:peano} and $\Ps(p)\succ0$, there exists $r>0$, smaller than the injectivity radius at $p$, such that $p$ is the unique minimizer of $\F$ on $K:=\overline B_r(p)$. Since $\F$ is continuous and $K$ is compact, for every $\varepsilon\in(0,r)$ the separation constant
\[
\alpha_\varepsilon:=
\inf_{\substack{x\in K\\\rho(x,p)\ge\varepsilon}}
\bigl(\F(x)-\F(p)\bigr)>0.
\]
This is precisely the well-separation condition in \cite[Theorem~5.7, p.~45]{vdV98}, restricted to $K$ and expressed for minimization. Next, writing $D=\operatorname{diam}(M)$, the triangle inequality yields $|\rho^2(x,y)-\rho^2(z,y)|\le 2D\rho(x,z)$ for all $x,z,y\in M$. Hence both $F_n$ and $\F$ are uniformly $2D$-Lipschitz. For any finite $\delta$-net $\{z_1,\ldots,z_N\}$ of $K$,
\[
\sup_{x\in K}|F_n(x)-\F(x)|
\le \max_{1\le j\le N}|F_n(z_j)-\F(z_j)|+4D\delta.
\]
The strong law of large numbers at the finitely many $z_j$, followed by $\delta\downarrow0$, gives $\eta_n:=\sup_{x\in K}|F_n(x)-\F(x)|\to0$ almost surely. Since $\hat p_n$ minimizes $F_n$ over $K$, we have $\F(\hat p_n)-\F(p)\le 2\eta_n$. Therefore $\rho(\hat p_n,p)\ge\varepsilon$ implies $\alpha_\varepsilon\le2\eta_n$, which is eventually impossible almost surely. Thus $\hat p_n\to p$ almost surely. This is the localized, almost-sure version of the argument in \cite[Theorem~5.7 and its proof, pp.~45--46]{vdV98}; that theorem itself states convergence in probability.

Consequently $\hat p_n$ eventually lies in the interior of $K$ and is a local empirical Fr\'echet minimizer. In the normal chart $\Theta=B_r(0)\subset T_pM$, $\hat\theta_n=\Log_p(\hat p_n)$ is consistent and, with probability tending to one, exactly maximizes the empirical criterion over $\Theta$. Thus the $o_P(n^{-1})$ near-maximization condition of \cite[Theorem~5.23, pp.~53--54]{vdV98} holds. The differentiability, Lipschitz, and quadratic-expansion conditions follow from parts (a)--(c) above, with $x_0$ replaced by $p$. Applying Theorem~5.23 gives the stated CLT.
\end{proof}

\vspace{-0.25em}

\noindent\emph{Scope.} The present argument avoids the \textit{separate} Donsker-class verification in Step~2 of the proof of \cite[Theorem~2, pp.~1240--1243]{HLW24} by applying \cite[Theorem~5.23, pp.~53--54]{vdV98}.\textit{ This does not eliminate empirical-process theory entirely}: the proof of Theorem~5.23 invokes \cite[Lemma~19.31, p.~284]{vdV98}, which establishes the stochastic approximation (19.30) using bracketing-entropy bounds from Example~19.7 and Theorem~19.28. It also invokes \cite[Corollary~5.53, p.~77]{vdV98}, whose proof uses bracketing-entropy estimates. The simplification lies in the verification of the hypotheses: the required $L^2(\mu)$-Lipschitz condition for the squared-distance criterion follows directly from compactness, whereas Step~2 of \cite[Theorem~2]{HLW24} establishes the Donsker property of the transported logarithm-field class using cut-locus-dependent estimates, particularly in Equations (44)--(47) there. The\textit{ local CLT} above concerns the specified \textit{local} empirical minimizers, not \textit{global }minimizers treated in \cite{HLW24}.
%%%%%%%%%%%%%%%%%%%%%%%%%%%%%%%%%%%%%%%%%%%%%%%%%%%%%%%%%%%%%%%%%%%%%%%%%%%%%%%%%%%%%%%%%%%%%%%%%%%%%%%%%%%%%%%%%%%%%%%%%%%%%%%%%%%%%%%%%%%%%%%%%%%%%%%%%%%%%%%%%%%%%%%%%%%%%%%%%%%%%%%%%%%%%%%%%%%%%%%%%%%

\section{Numerical verification of the cut locus corrected Hessian at a nonglobal local Fr\'echet minimum on $\mathbb{RP}^2$}
We illustrate Proposition~\ref{prop:stat} and the nondegeneracy condition in Corollary~\ref{cor:local} at a nonglobal local Fr\'echet minimum. Write $[y]=\{y,-y\}\in\mathbb{RP}^2$ for $y\in S^2$, so that $\rho([p],[y])=\arccos|p\cdot y|$. Let $Y\in S^2$ have density $f$ with respect to spherical area, where
\[
f(y)=\frac{\varepsilon}{4\pi}+\sum_{j=1}^3\frac{w_j\kappa_j}{4\pi\sinh\kappa_j}\cosh(\kappa_j a_j\cdot y).
\]
We take $\varepsilon=0.130$, $(w_1,w_2,w_3)=(0.405,0.405,0.060)$, $(\kappa_1,\kappa_2,\kappa_3)=(38,38,32)$, $a_1=e_3$, $a_2=e_1$, and $a_3=\sin(20^\circ)e_1+\cos(20^\circ)e_3$. The normalization follows from $\int_{S^2}\cosh(\kappa a\cdot y)\,d\sigma(y)=4\pi\sinh\kappa/\kappa$ for $\|a\|=1$ and $\varepsilon+\sum_jw_j=1$. Since $f(-y)=f(y)$, the pushforwarded/projected law $\mu$ of $\xi=[Y]$ has smooth density $\psi([y])=2f(y)$ with respect to Riemannian volume on $\mathbb{RP}^2$. In particular, this density is strictly positive on the cut locus of every point on $\mathbb{RP}^2.$

\begin{figure}[ht]
\centering
\includegraphics[width=\textwidth]{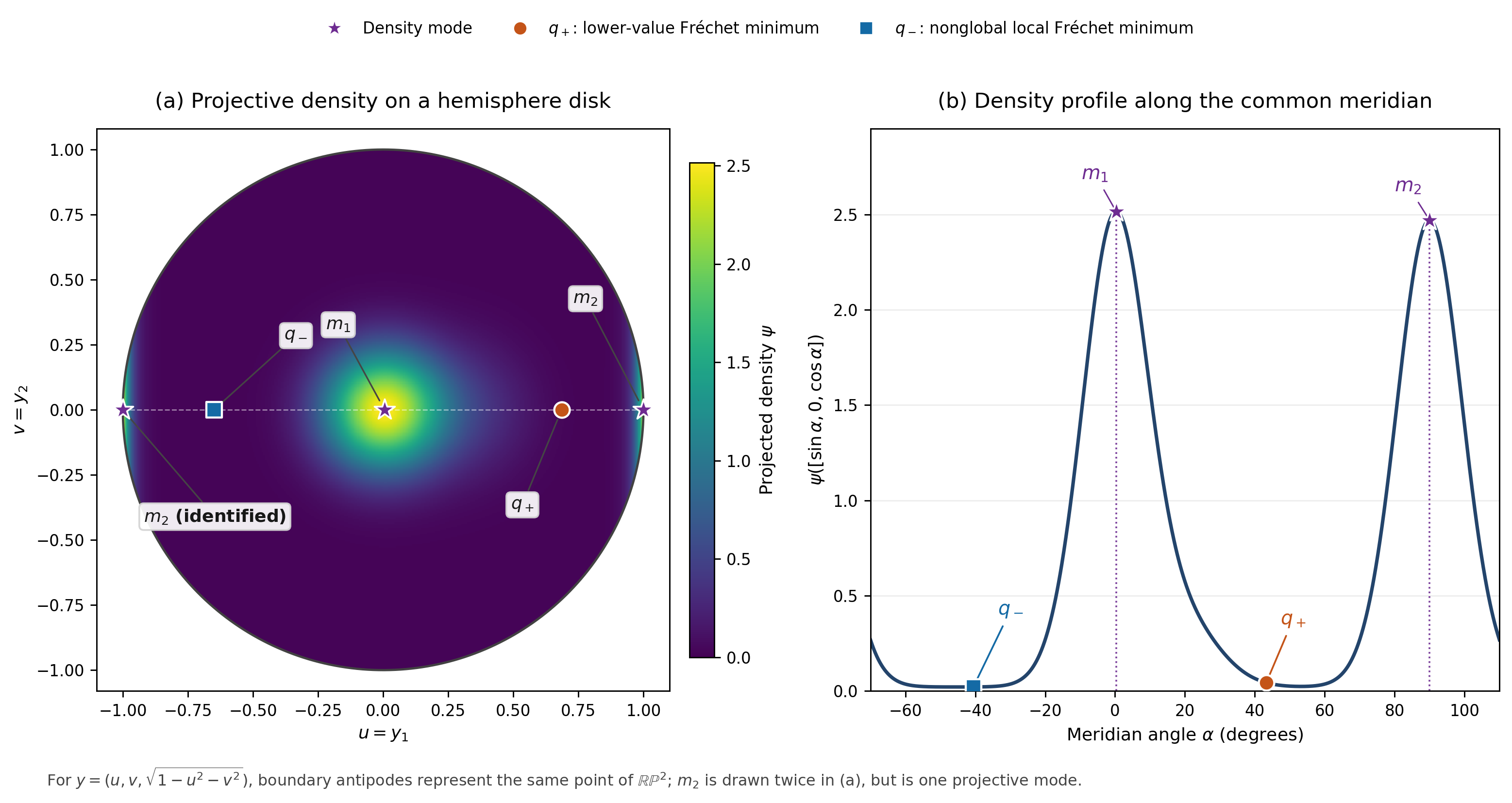}
\caption{Density on $\mathbb{RP}^2$ and its restriction to the great circle containing the two numerically identified Fr\'echet minima and the two modes.}
\label{fig:rp2-density}
\end{figure}

Besides antipodal symmetry, $f$ is invariant under the reflection $R(y_1,y_2,y_3)=(y_1,-y_2,y_3)$, since all three concentration directions $a_j$ lie in the $e_1e_3$-plane. Consequently, $\F([Rp])=\F([p])$, and both the density and the Fr\'echet function have vanishing transverse derivatives along the great circle $S^2\cap\{y_2=0\}$. The numerically identified density modes and the two Fr\'echet minima lie on this great circle. This symmetry explains why the derivatives perpendicular to the great circle vanish at these points, although it does not preclude additional extrema outside the circle. Figure~\ref{fig:rp2-density} illustrates the density on $\mathbb{RP}^2$, together with its restriction to the great circle containing the two numerically identified Fr\'echet minima.

Numerical minimization using the BFGS quasi-Newton algorithm of $\F$ from 60 random starting points identified two strict local minima near $q_\pm=[p_\pm]$, where $p_\pm=(\sin\alpha_\pm,0,\cos\alpha_\pm)$, $\alpha_+=43.274584^\circ$ and $\alpha_-=-40.548310^\circ$. The respective values $\F(q_+)=0.70072933$ and $\F(q_-)=0.76122945$, together with the numerically positive Hessians reported in the Table \ref{tab:rp2-local-hessian}, support that $q_-$ is a \textit{nonglobal local }minimum. The first point is the \textit{numerically identified global candidate}; global minimality or uniqueness is not claimed theoretically, only \textit{numerically}.

Let the \textit{expected pointwise Hessian of the squared distance} be $A_\mu(q)=\mathbb E[\operatorname{Hess}^M_q\rho^2(q,\xi)]$ and, for $\xi\notin C_q$, put the $\mathbb{RP}^2$-distance $r=\rho(q,\xi)$ and the unit vector $u=\Log_q\xi/r$. In an orthonormal tangent frame, $\operatorname{Hess}^M_q\rho^2(q,\xi)=2\{uu^{\mathsf T}+r\cot r(I-uu^{\mathsf T})\}$, with the removable value at $r=0$. For $q=[p]$, the explicit $\mathbb{RP}^2$ formula of \cite[\S2.4(c), p.~1230]{HLW24} becomes $J_\mu(q)=\pi\int_{S^2\cap p^\perp}f(y)\,yy^{\mathsf T}\,d\ell(y)$, where $yy^{\mathsf T}$ is an operator on $T_q\mathbb{RP}^2\simeq p^\perp$ and $d\ell$ is spherical arclength. Indeed, the integral over the full great-circle lift and the factor $\psi([y])=2f(y)$ reproduce their integral over the projective cut locus. 
Combining equations~(14)--(15) of \cite[pp.~1224--1225]{HLW24} with Proposition~\ref{prop:stat}, we numerically check the identity $\operatorname{Hess}^M\F(q)=A_\mu(q)-2J_\mu(q)$ at both identified \textit{local} minima.

\begin{table}[ht]
\centering\small\setlength{\tabcolsep}{4pt}
\begin{tabular}{@{}lcc@{}}\toprule
& $q_+$ (global candidate) & $q_-$ (nonglobal local) \\ \midrule
$A_\mu$ & $(1.92419748,\ 1.57564557)$ & $(1.92345299,\ 1.53076831)$ \\
$J_\mu$ & $(0.10224096,\ 0.10210694)$ & $(0.10955499,\ 0.10237353)$ \\
$A_\mu-2J_\mu$ & $(1.71971555,\ 1.37143170)$ & $(1.70434301,\ 1.32602126)$ \\
$\operatorname{Hess}^M\F$ (direct) & $(1.71971555,\ 1.37143170)$ & $(1.70434301,\ 1.32602126)$ \\ \bottomrule
\end{tabular}
\caption{Diagonal entries of the $2\times2$ operators in the orthonormal tangent frame $(\partial_\alpha p_\pm,e_2)$. The off-diagonal entries vanish by \textit{reflection symmetry }about the $e_1e_3$-plane. Numerical values are rounded.}\label{tab:rp2-local-hessian}
\end{table}

The correction $J_\mu$ was evaluated by numerical integration on the cut-locus, and also independently by an exact integration of the same formula given in \cite{HLW24}, pp. 1230, with agreement to floating-point precision. As an elementary check, the uniform component contributes exactly $\varepsilon\pi I/4$ to $J_\mu$. The expected pointwise Hessian $A_\mu$ was integrated separately, using 240 Gauss--Legendre nodes in $t=\cos r\in[0,1]$ and 512 angular nodes. Independently, $\operatorname{Hess}^M\F$ was computed from finite differences of the numerically integrated population function, using fourth-order centered differences for diagonal entries, centered differences for the mixed entry, and step $h=0.0025$. The unrounded maximum entrywise residuals in $\operatorname{Hess}^M\F-(A_\mu-2J_\mu)$ were $1.03\times10^{-10}$ at $q_+$ and $3.02\times10^{-10}$ at $q_-$. The corrected Hessians have numerically positive eigenvalues $(1.719716,1.371432)$ and $(1.704343,1.326021)$, respectively; hence $\Ps(q_\pm)=\tfrac12\operatorname{Hess}^M\F(q_\pm)$ is numerically positive definite. The experiment verifies the nonvanishing cut-locus contribution at a nonglobal strict local minimum to numerical accuracy, without claiming a proof of the local CLT or of global optimality.

\end{document}